\documentclass[12pt]{article}
\usepackage{amsmath, amssymb, amsthm, enumerate}
\usepackage[utf8]{inputenc}

\newtheorem{theorem}{Theorem}
\newtheorem{proposition}[theorem]{Proposition}

\newtheorem{lemma}[theorem]{Lemma}

\newtheorem{example}[theorem]{Example}
\newtheorem{question}[theorem]{Question}

\def\rr{{\mathbb R}}
\def\cc{{\mathbb C}}
\def\zz{{\mathbb Z}}
\def\qq{{\mathbb Q}}
\def\nnn{{\cal N}}
\def\pp{{\cal P}}
\def\al{\alpha}

\def\de{\delta}
\def\De{\Delta}

\def\si{\sigma}

\def\om{\omega}

\def\la{\lambda}

\def\su{\subset}

\def\cd{\cdot}
\def\stb{,\ldots ,}

\def\emp{\emptyset}

\def\ran{\rangle}
\def\lan{\langle}

\def\V{\Vert}

\def\sumis{\sum_{i=1}^s}

\def\sik{{\rr}^2}

\def\deg{{\rm deg}\, }

\def\msk{\medskip}
\def\bsk{\bigskip}
\def\noi{\noindent}

\begin{document}

\title{Separately polynomial functions}

\author{Gergely Kiss and Mikl\'os Laczkovich (Budapest, Hungary)}

\footnotetext[2]{{\bf MR subject classification:} 22A20, 54E45, 54E52}
\footnotetext[3]{Both authors were supported by the Hungarian National
Foundation for Scientific Research, Grant No. K124749. G. Kiss was supported by the Premium Postdoctoral Fellowship of the Hungarian Academy of Sciences.}

\maketitle

\begin{abstract} It is known that if $f\colon \sik \to \rr$ is a polynomial in
each variable, then $f$ is a polynomial. We present generalizations of this
fact, when $\sik$ is replaced by $G\times H$, where $G$ and $H$ are topological
Abelian groups. We show, e.g., that the conclusion holds (with generalized
polynomials in place of polynomials) if $G$ is a connected Baire space and $H$
has a dense subgroup of finite rank or, for continuous functions, if
$G$ and $H$ are connected Baire spaces. The condition of continuity can be
omitted if $G$ and $H$ are locally compact or one of them is metrizable. We
present several examples showing that the results are not far from being
optimal.
\end{abstract}

\section{Introduction}
It was proved by F. W. Carroll in \cite{C} that if $f\colon \sik \to \rr$ is a
polynomial in each variable, then $f$ is a polynomial. Our aim is to find
generalizations of this fact, when $\sik$ is replaced by the
product of two topological Abelian groups.

On topological Abelian groups we have to distinguish between the class of
polynomials and the wider class of generalized polynomials (see the next section
for the definitions). The two classes coincide if the group contains a dense
subgroup of finite rank. Now, the scalar product on the square of a Hilbert
space is an example of a continuous function which is a polynomial in each
variable, is a generalized polynomial on the product, but not a polynomial
(see Example \ref{ex1} below). Therefore, the appropriate problem is to find
conditions on the groups $G$ and $H$ ensuring that whenever a function on
$G\times H$ is a generalized polynomial in each variable, then it is a
generalized polynomial.

This problem was considered already by Mazur and Orlicz in \cite{MO} in the case
when $G$ and $H$ are topological vector spaces. They proved that if $X,Y,Z$ are
Banach spaces\footnote{Actually Mazur and Orlicz only assume that $X,Y,Z$
are $F$-spaces; that is, topological vector spaces whose topology is induced by
a complete invariant metric.} and the map $f\colon X\times Y\to Z$ is a
generalized polynomial in each variable, then $f$ is a generalized polynomial
\cite[Satz IV]{MO}. They also considered the case when continuity is not
assumed, and $X,Y,Z$ are linear spaces without topology \cite[Satz III]{MO}
(see also \cite[Lemma 1]{BS}). The topic has an extensive literature; see
\cite{Pl}, \cite{VKM} and the references therein.

In this note we consider the analogous problem when $G$ and $H$ are
topological Abelian groups. We show that if $G$ is a connected
Baire space, $H$ has a dense subgroup of finite rank, and if a function
$f\colon (G\times H) \to \cc$ is a generalized polynomial in each variable, then
$f$ is a generalized polynomial on $G\times H$ (Theorem \ref{t1}).
The same conclusion holds if $G$ and $H$ are both connected Baire spaces, and
one of them is metrizable or, if both are locally compact (Theorem \ref{t3}).

If $G$ and $H$ are connected Baire spaces, $f\colon (G\times H)\to \cc$ is
a generalized polynomial in each variable, and if $f$ has at least one
point of joint continuity, then $f$ is a generalized polynomial ((iii) of
Theorem \ref{t2}).

It is not clear if the extra condition of the existence of points of joint
continuity can be omitted from this statement (Question \ref{pr1}). The problem
is that a generalized polynomial must be continuous by definition, and a
separately continuous function on the product of Baire spaces can be
discontinuous everywhere, as it was shown recently in  \cite{MP}. In our case,
however, there are some extra conditions: the spaces are also connected, and the
function in question is a generalized polynomial. It is conceivable that
continuity follows under these conditions. As for the biadditive case, see
\cite{CF}.

There are several topological conditions implying that separately continuous
functions on a product must have points of joint continuity. In fact, the
topic has a vast literature starting with the paper \cite{N}. See, e.g., the
papers \cite{CM}, \cite{HW}, \cite{KM}, \cite{P}.

\section{Preliminaries}
Let $G$ be a topological Abelian group. We denote the group
operation by addition, and denote the unit by $0$. The translation operator
$T_h$ and the difference operator $\De _h$ are defined by $T_h f (x)= f(x+h)$
and $\De _h f(x)= f(x+h) -f(x)$ for every $f\colon G\to \cc$ and $h,x\in G$.

We say that a continuous function $f\colon G\to \cc$ is a {\it generalized
polynomial}, if there is an $n\ge 0$ such that $\De _{h_1} \ldots \De _{h_{n+1}} f
=0$ for every $h_1 \stb h_{n+1} \in G$. The smallest $n$ with this property is
the {\it degree of $f$}, denoted by $\deg f$. The degree of the identically zero
function is $-1$. We denote by ${\cal GP}={\cal GP}_G$ the set of generalized
polynomials defined on $G$.

A function $f\colon G\to \cc$ is said to be a {\it polynomial}, if there are
continuous additive functions $a_1 \stb a_n \colon G\to \cc$ and there is a
$P\in \cc [x_1 \stb x_n ]$ such that $f=P(a_1 \stb a_n )$. It is well-known that
every polynomial is a generalized polynomial. It is also easy
to see that the linear span of the translates of a polynomial is of finite
dimension. More precisely, a function is a polynomial if and only if it is a
generalized polynomial, and the linear span of its translates is of finite
dimension (see \cite[Proposition 5]{L}).  We denote by ${\cal P}={\cal P}_G$ the
set of polynomials defined on $G$.

Let $f$ be a complex valued function defined on $X\times Y$. The {\it sections}
$f_x \colon Y\to \cc$ and $f^y \colon X\to \cc$ of $f$ are defined by
$f_x (y)=f^y (x)=f(x,y)$ $(x\in X, \ y\in Y)$.

Let $G,H$ be topological Abelian groups. A function $f\colon
(G\times H) \to \cc$ is a {\it separately polynomial function} if $f_x \in
\pp _H$ for every $x\in G$ and $f^y \in \pp _G$ for every $y\in H$. Similarly,
we say that $f\colon (G\times H) \to \cc$ is a {\it separately generalized
polynomial function} if $f_x \in {\cal GP} _H$ for every $x\in G$ and $f^y
\in {\cal GP}_G$ for every $y\in H$.

In general we cannot expect that every separately polynomial function on
$G\times H$ is a polynomial; not even if $G=H$ is a Hilbert space.

\begin{example} \label{ex1}
{\rm Let $G$ be the additive group of an infinite dimensional Hilbert space.
Then the scalar product $f(x,y)=\lan x,y\ran$ on $G^2$ is a separately
polynomial function, since its sections are continuous additive functions. In
fact, $f^y$ is a linear functional and $f_x$ is a conjugate linear functional
for every $x,y\in G$. Thus the sections of $f$ are polynomials.

Now, while  the scalar product is a generalized polynomial (of degree $2$) on
$G^2$, it is not a polynomial on $G^2$, because the dimension of the linear
span of its translates is infinite. Indeed, let $g(x)=\lan x,x\ran =\V x\V ^2$
for every $x\in G$. Then $\De _h g(x)=2\lan h,x\ran +\V h\V ^2$ for every
$h\in G$. It is easy to see that the functions $\lan h,x\ran$ $(h\in G)$
generate a linear space of infinite dimension, and then the same is true for
the translates of $g$ and then for those of $f$ as well.
}
\end{example}

Therefore, the best we can expect is that, under suitable conditions on $G$ and
$H$, every separately generalized polynomial function on $G\times H$ is a
generalized polynomial.

We denote by $r_0 (G)$ the torsion free rank of the group $G$; that is, the
cardinality of a maximal independent system of elements of $G$ of infinite
order. Thus $r_0 (G)=0$ if and only if $G$ is torsion. In the sequel by the rank
of a group we shall mean the torsion free rank.
It is known that if $G$ has a dense subgroup of finite rank, then the classes
of polynomials and of generalized polynomials on $G$ coincide (see
\cite[Theorem 9]{L}).

The set of roots of a function $f\colon G\to \cc$ is denoted by $Z_f$. That is,
$Z_f =\{ x\in G \colon f(x)=0\}$. We put
$${\cal N}_P = {\cal N}_P (G) =\{ A\subset G\colon \exists p\in {\cal P}_G ,
\ p\ne 0 , \ A\su Z_p \}$$
and
$${\cal N}_{GP} = {\cal N}_{GP} (G) = \{ A\subset G\colon \exists p\in
{\cal GP}_G ,\ p\ne 0, \ A\su Z_p \} .$$
It is easy to see that ${\cal N}_P$ and ${\cal N}_{GP}$ are proper ideals of
subsets of $G$. Let $\nnn_P^\si$ and $\nnn_{GP}^\si$ denote the $\si$-ideals
generated by $\nnn _P$ and $\nnn _{GP}$, respectively. Note that $\nnn _P \su
\nnn_{GP}$ and $\nnn _P^\si  \su \nnn_{GP}^\si$.

If $G$ is discrete, then $\nnn_P^\si$ and $\nnn_{GP}^\si$ are not proper
$\si$-ideals (except when $G$ is torsion), according to the next observation.
\begin{proposition}\label{p1}
Let $G$ be a discrete Abelian group. If $G$ is not torsion, then $G\in
\nnn_P^\si$.
\end{proposition}
\noindent
{\bf Proof.} Let $a\in G$ be an element of infinite order. Then $\phi (na)=n$
$(n\in \zz )$ defines a homomorphism from the subgroup generated by $a$ into
$\qq$, the additive group of the rationals. Since $\qq$ is divisible, $\phi$ can
be extended to $G$ as a homomorphism from $G$ into $\qq$. Let $\psi$ be such an
extension.

Then $p_r =\psi +r$ is a nonzero polynomial on $G$ for every $r\in \qq$.
If $x\in G$, then $x$ is the root of $p_r$, where $r=-\psi (x)\in \qq$.
Therefore, $G=\bigcup_{r\in \qq} Z_{p_r} \in \nnn_P^\si$. \hfill $\square$

\msk
A simple sufficient condition for $G\notin \nnn_{GP}^\si$ is given by the
next result.
\begin{lemma}\label{l1}
If $G$ is a connected Baire space, then the $\si$-ideals $\nnn_P^\si$  and $\nnn_{GP}^\si$ are proper; that is, $G\notin \nnn_P^\si$ and $G\notin \nnn_{GP}^\si$.
\end{lemma}

\noindent
{\bf Proof.}
It is enough to prove that every element of $\nnn_{GP}$ is nowhere dense.
Suppose $A\in \nnn_{GP}$ is dense in a nonempty open set $U$. Let $p\in
{\cal GP} (G)$ be a nonzero generalized polynomial vanishing on $A$. Since
$A\su Z_p$ and $Z_p$ is closed, we have $U\su Z_p$. Since $G$ is connected,
every neighbourhood of the origin generates $G$. It is known that in such a
group, if a generalized polynomial vanishes on a nonempty open set, then it
vanishes everywhere (see \cite[Theorem 3.2, p.~33]{Sz}).
This implies that $p$ is identically zero, which is impossible. \hfill $\square$

\section{Main results}

\begin{theorem} \label{t1}
Let $G,H$ be topological Abelian groups, and suppose that
\begin{enumerate}[{\rm (i)}]
\item $\nnn ^\si _{GP} (G)$ is a proper $\sigma$-ideal in $G$, and
\item $H$ has a dense subgroup of finite rank.
\end{enumerate}
If $f\colon (G\times H)\to \cc$ is a
separately generalized polynomial function, then $f$ is a generalized
polynomial on $G\times H$.
\end{theorem}

\begin{theorem} \label{t2}
Let $G,H$ be topological Abelian groups, and suppose that $\nnn ^\si _{GP} (G)$
is a proper $\sigma$-ideal in $G$, and $\nnn ^\si _{GP} (H)$ is a proper
$\sigma$-ideal in $H$. Then the following statements are true.
\begin{enumerate}[{\rm (i)}]
\item If $f\colon (G\times H) \to \cc$ is a separately generalized polynomial
function, then $f$ is a generalized polynomial on $G\times H$ with respect to
the discrete topology.
\item Every joint continuous separately generalized polynomial function
$f\colon (G\times H) \to \cc$ is a generalized polynomial on $G\times H$.
\item If $G$ and $H$ are connected and a separately generalized polynomial
function $f\colon (G\times H) \to \cc$ has at least one point of joint
continuity, then $f$ is a generalized polynomial on $G\times H$.
\end{enumerate}
\end{theorem}

By Lemma \ref{l1}, (i) of Theorem \ref{t1} can be replaced by the
condition that $G$ is a connected Baire space. Similarly, the condition of
Theorem \ref{t2} can be replaced by the condition that $G$ and $H$ are
connected Baire spaces.

As for (iii) of Theorem \ref{t2} note the following facts.
\begin{enumerate}[$\bullet$]
\item If $X,Y$ are nonempty topological spaces, $X$ is Baire, $Y$ is first
countable and $f\colon X\times Y\to \cc$ is
separately continuous, then $f$ has at least one point of joint continuity.
(See, e.g. \cite[p.~441]{W}.)
\item A topological group is first countable if and only if it is metrizable.
\item If $X,Y$ are nonempty locally compact and $\si$-compact topological
spaces, $f\colon X\times Y\to \cc$ is separately continuous, then $f$ has at least one point of joint continuity. (See \cite[Theorem 1.2]{N}).
\item Every connected and locally compact topological group is $\si$-compact.
\end{enumerate}

Comparing these with (iii) of Theorem \ref{t2} we obtain the following.

\begin{theorem} \label{t3}
Suppose that the topological Abelian groups $G,H$ are connected and Baire,
and either
\begin{enumerate}[{\rm (i)}]
\item at least one of $G$ and $H$ is metrizable, or
\item $G$ and $H$ are locally compact.
\end{enumerate}
If $f\colon (G\times H) \to \cc$ is a separately generalized polynomial
function, then $f$ is a generalized polynomial on $G\times H$. \hfill $\square$
\end{theorem}

\begin{question} \label{pr1}
{\rm Are the conditions (i) and (ii) necessary in the statement of
Theorem \ref{t3}? (See the introduction.)}
\end{question}

We prove Theorems \ref{t1} and \ref{t2} in the next section. In Section 5
we present examples showing that some of the conditions appearing in
Theorems \ref{t1} and \ref{t2} cannot be omitted.

\section{Proof of Theorems \ref{t1} and \ref{t2}}

\begin{lemma}\label{l2}
Let $H$ be a topological Abelian group, and suppose that $H$ has a dense
subgroup of finite rank. Then, for every positive integer $d$, there are
finitely many points $x_1 \stb x_s \in H$ and there are generalized polynomials
$q_1 \stb q_s \in {\cal GP}_H$ of degree $<d$ such that $p=\sum_{i=1}^s p(x_i )
\cd q_i$ for every $p\in {\cal GP}_H$ with $\deg p < d$.
\end{lemma}

\noindent
{\bf Proof.} Let ${\cal GP}^{<d}$ denote the set of generalized polynomials
$f\in {\cal GP}_H$ of degree $<d$. Clearly, ${\cal GP}^{<d}$ is a linear space
over $\cc$.

Let $K$ be a dense subgroup of $H$ with $r_0 (K) =N<\infty$.
Let $\{ h_1 \stb h_N \}$ be a maximal set of independent elements of $K$ of
infinite order, and let $L$ denote the subgroup of $K$
generated by the elements $h_1 \stb h_N$. If $k=(k_1 \stb k_N )\in \zz ^N$, then
we put $\V k\V =\max_{1\le i\le N} |k_i |$.
We abbreviate the sum $\sum_{i=1}^N k_i \cd h_i$ by $\lan k, h \ran$. Then we
have $L=\{ \lan k, h \ran \colon k\in \zz ^N\}$. We put
$$A=\{ \lan k, h\ran \colon k\in \zz ^N , \ \V k \V \le [d/2] \} .$$
First we prove that if $p\in {\cal GP}^{<d}$ vanishes on $A$, then $p=0$.

Suppose $p\ne 0$. Since $p$ is continuous and $K$ is dense in $H$, there is an
$x_0 \in K$ such that $p(x_0 )\ne 0$. The maximality of the system
$\{ h_1 \stb h_N \}$ implies that $nx_0 \in L$ with a suitable nonzero integer
$n$. It is easy to see that there is a polynomial $P\in \cc [x]$ such that
$p(mx_0 )=P(m)$ for every integer $m$. Since $P(1)=p(x_0 )\ne 0$, it follows
that $P\ne 0$, hence $P$ only has a finite number of roots. Thus $p(mnx_0 )=
P(mn)\ne 0$ for all but a finite number of integers $m$. Fix such an $m$. Then
$mnx_0 \in L$, and thus $mnx_0  =\lan k, h \ran$ with a suitable $k\in \zz ^N$.
We find that $p(\lan k, h \ran )\ne 0$ for some $k\in \zz ^N$.

Let $k=(k_1 \stb k_N )\in \zz ^N$ be such that $p(\lan k, h \ran )\ne 0$ and
$\V k\V$ is minimal. If $\V k\V \le [d/2]$, then $\lan k, h \ran \in A$, and we
have $p(\lan k, h \ran )=0$ by assumption. Thus we have $\V k \V >[d/2]$.
Put $\ell =(\ell _1 \stb \ell_N )$, where
$$\ell _i =
\begin{cases}
1 & \text{if $k_i > [d/2]$,} \cr
0 & \text{if $|k_i |\le [d/2]$,} \cr
-1 & \text{if $k_i <-[d/2]$} \cr
\end{cases}
\qquad (i=1\stb N).$$
Then we have $\V k -j \ell \V <\V k\V$ for every $j=1\stb d$.
By the minimality of $\V k\V$ we have
$p( \lan k-j \ell ,h \ran )=0$ for every $j=1\stb d$.

Put $v=\lan \ell , h \ran$. Since $\deg p< d$, it follows that $\De _{-v}^{d}
p(x) =0$ for every $x\in H$. Now we have
\begin{align*}
0&=\De _{-v}^{d} p (\lan k, h \ran )=\sum_{j=0}^{d} (-1)^{d-j} {d \choose j}
p(\lan k, h \ran -j v) =\\
&= (-1)^{d} p(\lan k, h \ran )+\sum_{j=1}^{d} (-1)^{d-j}
{d \choose j} p(\lan k-j\ell , h \ran )=\\
&=(-1)^{d} p(\lan k, h \ran ),
\end{align*}
which is impossible. This proves $p=0$.

The set of functions $V=\{ p|_A \colon p\in {\cal GP}^{<d} \}$ is a
finite dimensional linear space
over $\cc$. The map $p\mapsto p|_A$ is linear from ${\cal GP}^{<d}$ onto $V$
and, as we proved above, it is injective. Therefore, ${\cal GP}^{<d}$ is of
finite dimension.

Let $b_1 \stb b_s$ be a basis of ${\cal GP}^{<d}$. Since the functions
$b_1 \stb b_s$ are linearly independent, there are elements $x_1 \stb x_s$ such
that the determinant $\det |b_i (x_j )|$ is non\-zero (see \cite[Lemma 1,
p. 229]{AD}). Put $X=\{ x_1 \stb x_s \}$. Then $b_1 |_X \stb b_s |_X$ are
linearly independent, and thus the map $f\mapsto f|_X$ is bijective and linear
from ${\cal GP}^{<d}$ onto the set of functions $f\colon X\to \cc$.

Then there are functions $q_1 \stb q_s \in {\cal GP}^{<d}$ such that $q_i (x_i )=
1$ and $q_i (x_j )=0$ for every $i,j=1\stb s, \ i\ne j$.

Let $p\in {\cal GP}^{<d}$ be given. Then $p-\sumis p(x_j ) q_j$ is a generalized
polynomial of degree $<d$ vanishing on $X$, hence on $H$. That is, we have
$p=\sumis p(x_j ) q_j$. \hfill $\square$

\msk \noi
{\bf Proof of Theorem \ref{t1}.}
Let $f\colon (G\times H) \to \cc$ be a separately generalized polynomial
function. Put $G_n =\{ x\in G \colon \deg f_x <n \}$ $(n=1,2, \ldots )$.
Since $\nnn ^\si _{GP} (G)$ is a proper $\sigma$-ideal in $G$, there is an $n$
such that $G_n \notin \nnn _{GP} (G)$. Fix such an $n$.

By Lemma \ref{l2}, there are points $y_1 \stb y_s \in H$ and generalized
polynomials $q_1 \stb q_s \in {\cal GP}_H$ such that $p=\sum_{i=1}^s p(y_i )
\cd q_i$ for every $p\in {\cal GP}_H$ with $\deg p < n$. Therefore, we have
$$f(x,y)=\sumis f(x,y_i )q_i (y)$$
for every $x\in G_n$ and $y\in H$. If $y\in H$ is fixed, then
$f(x,y)-\sumis f(x,y_i )q_i (y)$ is a generalized polynomial on $G$ vanishing
on $G_n$. Since $G_n \notin \nnn _{GP} (G)$, it follows that
$f(x,y)=\sumis f(x,y_i )q_i (y)$ for every $(x,y)\in G\times H$.
By $f^{y_i} \in {\cal GP}_G$ and $q_i \in {\cal GP}_H$, we obtain $f\in
{\cal GP}_{G\times H}$. \hfill $\square$

\bsk \noindent
{\bf Proof of Theorem \ref{t2}.}

\noi
(i) Suppose $f$ satisfies the conditions. By Lemma \ref{l3}, it is enough to
show that the degrees $\deg f_x$ and $f^y$ are bounded.

Put $A_n =\{ x\in G\colon \deg f_x < n\}$. Then $G=\bigcup_{n=1}^\infty A_n$.
Since $\nnn ^\si_{GP}(G)$ is a proper $\sigma$-ideal, there is an $n$ such that
$A_n \notin \nnn _{GP} (G)$. We fix such an $n$, and prove that
\begin{equation} \label{e5}
\De _{(0,h_1 )} \ldots \De _{(0,h_n )} f=0
\end{equation}
for every $h_1 \stb h_n \in H$.
Let $g$ denote the left hand side of
\eqref{e5}. Then $g(x,y)=\sum_{i=1}^s a_i f(x,y+b_i )$, where $s=2^n$,
$a_i =\pm 1$ and $b_i \in H$ for every $i$. Let $y\in H$ be fixed. Then
$g^y = \sum_{i=1}^s a_i f^{y+b_i}$, and thus $g^y$ is a generalized polynomial
on $G$.

If $x\in A_n$, then $\deg f_x <n$, and thus $g_x =0$. Therefore $g^y (x)=0$
for every $x\in A_n$. Since $g^y$ is a generalized polynomial and
$A_n \notin \nnn _{GP}(G)$, it follows that $g^y =0$. Since $y$ was arbitrary, this
proves \eqref{e5}. Thus $\deg f_x <n$ for every $x\in G$.

A similar argument shows that, for a suitable $m$, $\deg f^y <m$ for every
$y\in H$.

\noi
Statement (ii) of the theorem is clear from (i).

Suppose that $G$ and $H$ are connected. Now we use the fact
that if $f$ is a discrete generalized polynomial on an Abelian group which is
generated by every neighbourhood of the origin, and if $f$ has a point of
continuity, then $f$ is continuous everywhere. (See \cite[Theorem 3.6]{Sz}) or,
for topological vector spaces, \cite[Theorem 1]{BS}.)
In our case the group $G\times H$ is connected, so the condition is satisfied,
and we conclude that $f$ is continuous everywhere on $G\times H$.
Thus (iii) follows from (ii). \hfill $\square$

\section{Examples}
First we show that in Theorem \ref{t1} none of the conditions on $G$ and $H$
can be omitted. First we show that without condition (i) the conclusion of
Theorem \ref{t1} may fail. We shall need the easy direction of the following
result.
\begin{lemma}\label{l3}
Let $G,H$ be discrete Abelian groups. A function $f\colon (G\times H) \to \cc$
is a generalized polynomial if and only if the sections $f_x$ $(x\in G)$ and
$f^y$ $(y\in H)$ are generalized polynomials of bounded degree.
\end{lemma}

\noindent
{\bf Proof.}
Suppose $f\colon (G\times H) \to \cc$ is a generalized polynomial of degree
$<d$. Then $\De _{(x_1 ,0)} \ldots \De _{(x_d ,0)} f=0$ for every $x_1 \stb x_d \in
G$. Then, for every $y\in H$, we have $\De _{x_1} \ldots \De _{x_d} f^y=0$
for every $x_1 \stb x_d \in G$, and thus $f^y$ is a generalized polynomial of
degree $<d$ for every $y\in H$. A similar argument shows that $f_x$ is a
generalized polynomial of degree $<d$ for every $x\in G$, proving the
``only if'' statement.

Now suppose that $f\colon (G\times H) \to \cc$ is such that
$f_x$ $(x\in G)$ and $f^y$ $(y\in H)$ are generalized polynomials of
degree $<d$. Then we have
\begin{equation} \label{e1}
\De _{(h_1 ,0)} \ldots \De _{(h_d ,0)} f=0
\end{equation}
for every $h_1 \stb h_d \in G$, and
\begin{equation} \label{e2}
\De _{(0,k_1 )} \ldots \De _{(0,k_d )} f=0
\end{equation}
for every $k_1 \stb k_d \in H$. In order to prove that $f$ is a
generalized polynomial of  degree $<2d$, it is enough to show that
\begin{equation} \label{e3}
\De _{(a_1 ,b_1 )} \ldots \De _{(a_{2d} ,b_{2d})} f=0
\end{equation}
for every $(a_i ,b_i )\in G\times H$ $(i=1\stb 2d)$.
The identity $\De_{u+v}=T_u \De_v  +\De_u$ gives
$$\De _{(a_i ,b_i )} =T_{(a_i ,0)} \De _{(0,b_i )} +\De _{(a_i ,0)}$$
for every $i$. Therefore, the left hand side of \eqref{e3} is the sum of terms
of the form $T_c \De_{c_1} \ldots \De_{c_{2d}} f$, where $c\in G\times \{ 0\}$, and
$c_i \in (G\times \{ 0\}) \cup (\{ 0\} \times H)$ for every $i$.
If there are at least $d$ indices $i$ with $c_i \in (G\times \{ 0\})$, then
\eqref{e1} gives $ \De_{c_1} \ldots \De_{c_{2d}} f=0$. Otherwise there
are at least $d$ indices $i$ with $c_i \in (\{ 0\} \times H)$, and then
\eqref{e2} gives $ \De_{c_1} \ldots \De_{c_{2d}} f=0$. This proves
\eqref{e3}. \hfill $\square$

\msk
Now we turn to the first example.
\begin{example} \label{ex2}
{\rm Let $G,H$ be discrete Abelian groups. We show that if none of $G$ and $H$
is torsion, then there is a separately polynomial function
$f\colon (G\times H) \to \cc$ such that $f$ is not a generalized
polynomial on $G\times H$.

By Proposition \ref{p1}, $\nnn^\si_{P} (G)$ is not a proper $\sigma$-ideal; that
is, $G= \bigcup_{n=1}^\infty A_n$, where $A_n \ne \emp$ and $A_n \in \nnn _{P} (G)$
for every $n$. Let $p_n \in {\cal P}_G$ be such that $p_n \ne 0$ and $A_n \su
Z_{p_n}$. Then $p_n$ is not constant; that is, $\deg p_n \ge 1$.

Let $P_n =p_1 \cdots p_n$; then $P_n (x)=0$ for every $x\in \bigcup_{i=1}^n A_i$,
and we have $0<\deg P_1 <\deg P_2 < \ldots$. (Here we use the fact that
$\deg pq=\deg p +\deg q$ for every $p,q\in {\cal GP}_G$, $p,q\ne 0$.)
Note that for every $x\in G$ we have $P_n (x)=0$ for all but a finite number of
indices $n$.

Similarly, we find polynomials $Q_n \in \pp _H$ such that
$0<\deg Q_1 <\deg Q_2 < \ldots$, and for every $y\in H$ we have $Q_n (y)=0$ for
all but a finite number of indices $n$.

We put $f(x,y)=\sum _{n=1}^\infty P_n (x)Q_n (y)$ for every $x\in G$ and $y\in H$.
If $y\in H$ is fixed, then the sum defining $f$ is finite, and thus $f^y
\in {\cal P}_G$. Similarly, we have $f_x \in {\cal P}_H$ for every $x\in G$.

The degrees $\deg f^y$ $(y\in H)$ are not bounded. Indeed, for every $N$, there
is an $y\in H$ such that $Q_N (y)\ne 0$. Then
$f^y = \sum _{n=1}^M  Q_n (y) \cd P_n$ with an $M\ge N$, where the coefficients
$Q_n (y)$ are nonzero if $n\le N$. Therefore, $\deg f^y \ge \deg P_N \ge N$,
proving that the set $\{ \deg f^y \colon y\in H\}$ is not bounded. By Lemma
\ref{l3}, it follows that $f$ is a not a generalized polynomial.
}
\end{example}

By the example above, if $G$ and $H$ are discrete Abelian groups of
positive and finite rank, then the conclusion of Theorem \ref{t1} fails.
That is, $G\notin \nnn ^\si _{GP} (G)$ cannot be omitted from the conditions
of Theorem \ref{t1}.

Next we show that the condition on $H$ cannot be omitted either.

\begin{example} \label{ex3}
{\rm Let $H$ be a discrete Abelian group of infinite rank. We show that
if $G$ is a topological Abelian group such that $\pp _G$ contains nonconstant
polynomials, then there is a continuous separately polynomial function $f$ on
$G\times H$ such that $f$ is not a generalized polynomial.

Let $h_\al$ $(\al <\kappa )$ be a maximal set of independent elements of $H$ of
infinite order, where $\kappa \ge \om$. Let $K$ denote the subgroup of $H$
generated by the elements $h_\al$ $(\al <\kappa )$. Every element of $K$
is of the form $\sum_{\al <\kappa}  k_\al h_\al$, where $k_\al \in \zz$
for every $\al$, and all but a finite number of the coefficients $k_\al$ equal
zero.

Let $p \in \pp _G$ be a nonconstant polynomial. We define
$f(x,y)=\sum_{i=1}^\infty k_i \cd p^i (x)$ for every $x\in G$ and $y\in K$,
$y=\sum_{\al <\kappa}  k_\al h_\al$. (Note that the sum only contains a finite
number of nonzero terms for every $x$ and $y$.) In this way we defined $f$ on
$G\times K$ such that  $f_x$ is additive on $K$ for every $x\in G$.

If $y\in H$, then there is a nonzero integer $n$ such that $ny\in K$. Then we
define $f(x,y)=\tfrac{1}{n}\cd f(x,ny)$ for every $x\in G$. It is easy to see
that $f(x,y)$ is well-defined on $G\times H$, and $f_x$ is additive on $H$ for
every $x\in G$. Therefore, $f_x$ is a polynomial on $G$ for every $x\in G$.

If $y\in H$ and $ny\in K$ for a nonzero integer $n$, then $f^y$ is of the
form $\tfrac{1}{n}\cd \sum_{i=1}^N k_i \cd p^i$, and thus $f^y \in \pp _G$.
Since $f^y$ is continuous for every $y\in H$ and $H$ is discrete, it follows
that $f$ is continuous on $G\times H$.

Still, $f$ is not a generalized polynomial on $G\times H$, as the set of degrees
$\deg f^y$ $(y\in H)$ is not bounded: if $y=h_i$, then
$f^y =p^i$, and $\deg p^i =i\cd \deg p \ge i$ for every $(i=1,2,\ldots )$.
}
\end{example}

In the example above we may choose $G$ in such a way that $G\notin
\nnn ^\si _{GP} (G)$ holds. (Take, e.g., $G=\rr$.) In our next example this
condition holds for both $G$ and $H$.

\begin{example} \label{ex4}
{\rm Let $E$ be a Banach space of infinite dimension, and let $G$ be the
additive group of $E$ equipped with the weak topology $\tau$ of $E$.
It is well-known that every ball in $E$
is nowhere dense w.r.t.~$\tau$, and thus $G$ is of first category in itself.

Still, we show that $G\notin \nnn ^\si _{GP} (G)$. Indeed, the original norm
topology of $E$ is stronger than $\tau$, and makes $E$ a connected Baire space.
If a function is continuous w.r.t.~$\tau$, then it is also continuous
w.r.t.~the norm topology. Therefore, every polynomial $p\in {\cal P}(G)$ is
also a polynomial on $E$, and thus $\nnn _P (G)\su \nnn _P (E)$ and
$\nnn ^\si_P (G)\su \nnn ^\si_P (E)$. Since $\nnn ^\si_P (E)$ is proper by
Lemma \ref{l1}, it follows that $\nnn ^\si_P (G)$ is proper. The same is true
for $\nnn ^\si_{GP} (G)$.

Now let $H$ be an infinite dimensional Hilbert space, and let $G$ be the
additive group of $H$ equipped with the weak topology of $H$.
Let $f$ be the scalar product on $H^2$. Since the linear functionals and
conjugate linear functionals are continuous w.r.t.~the weak topology,
it follows that $f$ is a separately polynomial function on $G^2$ (see Example
\ref{ex1}).

However, $f$ is not a generalized polynomial on $G^2$, since $f$ is not
continuous. In order to prove this, it is enough to show that $f(x,x)=
\V x\V ^2$ is not continuous on $H$ w.r.t.~the weak topology. Suppose it is.
Then there is a neighbourhood $U$ of $0$ such that $\V x\V <1$ for every
$x\in U$. By the definition of the weak topology, there are linear
functionals $L_1 \stb L_n$ and there is a $\de >0$ such that whenever
$|L_i (x)|<\de$ $(i=1\stb n)$, then $\V x\V <1$.

Since $H$ is of infinite dimension, there is an $x\ne 0$ such that
$L_i (x)=0$ for every $i=1\stb n$. (Otherwise every linear functional would
be a linear combination of $L_1 \stb L_n$, and then $H=H^*$ would be finite
dimensional.) Then $\la x\in U$ for every $\la \in \cc$
and $\V \la x\V <1$ for every $\la \in \cc$, which is impossible.
}
\end{example}

The example above shows that in (ii) of Theorem \ref{t2} the condition
of joint continuity cannot be omitted. Note also that the group $G$ defined in
Example \ref{ex4} is a topological vector space, hence connected. This shows
that in (iii) of Theorem \ref{t2} the condition of the existence of points of
joint continuity cannot be omitted either.

\noi
{\bf Acknowledgment.} We are indebted to the referee for calling our
attention to important pieces of literature and for several suggestions
that improved the paper considerably.

\end{document}